\documentclass[11pt]{amsart}

\usepackage{amsthm,amsmath,amssymb,amscd,amsfonts,latexsym}

\theoremstyle{plain}
\newtheorem{theorem}{Theorem}[section]

\newtheorem{proposition}[theorem]{Proposition}
\newtheorem{corollary}[theorem]{Corollary}
\newtheorem{lemma}[theorem]{Lemma}

\theoremstyle{definition}
\newtheorem{definition}[theorem]{Definition}

\newtheorem{remark}[theorem]{Remark}

\newtheorem*{acknowledgement}{Acknowledgement}

\newcommand{\PP}{\mathbb{P}}

\newcommand{\NN}{\mathbb{N}}

\newcommand{\CC}{\mathbb{C}}

\newcommand{\calo}{{\mathcal O}}
\newcommand{\OO}{{\mathcal O}}
\newcommand{\TT}{{\mathcal T}}

\newcommand{\frakm}{{\mathfrak m}}
\newcommand{\fraka}{{\mathfrak a}}
\newcommand{\fra}{{\mathfrak a}}
\newcommand{\frakb}{{\mathfrak b}}

\newcommand{\Bl}{\text{Bl}}
\newcommand{\dra}{\dashrightarrow}
\newcommand{\lra}{\longrightarrow}

\DeclareMathOperator{\ord}{ord}

\DeclareMathOperator{\divisor}{div}
\DeclareMathOperator{\Dom}{Dom}
\DeclareMathOperator{\dbar}{{\overline \partial}}
\DeclareMathOperator{\rp}{Re}
\DeclareMathOperator{\Zeroes}{Zeroes}

\begin{document}

\title{Finite type and the effective Nullstellensatz}

\begin{abstract} 
Improved local and global versions of the effective Nullstellensatz for ideal sheaves on non-singular complex varieties are obtained, based on a new invariant motivated by the notion of finite type from the theory of several complex variables. Two closely related curve selection theorems for curves with maximal adjusted tangency orders to a given ideal sheaf are established along the way, using normalized blow-ups and integral closures.  
\end{abstract}

\author[G. Heier]{Gordon Heier}
\address{Ruhr-Universit\"at Bochum\\
Fakult\"at f\"ur Mathematik\\
44780 Bochum\\
Germany}
\curraddr{University of Michigan, Mathematics Department, 2074 East Hall, 530 Church Street, Ann Arbor, MI 48109, USA}
\email{heier@math.harvard.edu}

\thanks{This research was supported by the DFG Schwerpunkt {\it Global methods in complex geometry}.}
 

\maketitle

\section{Introduction} 
The notion of {\it finite type} plays a crucial role in the $\dbar$-Neumann problem in several complex variables. In this Introduction, we will provide a short survey aimed at readers unfamiliar with finite type, subelliptic estimates, and subelliptic multiplier ideal sheaves. The subsequent sections are not logically dependent on the analytic phenomena described. Section \ref{curve_selection} contains the statement and proof of the curve selection theorem for finite type ideals from \cite{HL}, which has become part of this article (see page \pageref{ackn}). Section \ref{eff_nsts} makes the connection with the Nullstellensatz in the spirit of \cite{Ein_Lazarsfeld}. The new angle is that, based on the notion of finite type, a new invariant for an ideal sheaf is defined that is well suited for the effective Nullstellensatz.\par
\subsection{Finite type domains and subelliptic estimates}
Let $\Omega\subset \CC^n$ be a weakly pseudoconvex domain. We assume that $\Omega$ is given as $\{r<0\}$ for some real-valued $C^\infty$ function $r$ whose gradient does not vanish anywhere on $\{r=0\}$.\par
Kohn and Nirenberg (\cite{Kohn_Nirenberg}) showed that local regularity for the $\dbar$-Neumann problem at any point $p\in \overline \Omega$  is a consequence of a subelliptic a priori estimate, defined as follows. As usual, we restrict ourselves to the case of $(0,1)$-forms. 
\begin{definition}
The $\dbar$-Neumann problem for $(0,1)$-forms is said to satisfy a {\it subelliptic estimate (of order $\varepsilon$)} at $p$ if there exists a neighborhood $U$ of $p$, a constant $\varepsilon > 0$ and a constant $C>0$ such that 
\begin{equation*}
\|\varphi\|^2_\varepsilon \leq C(\|\dbar \varphi\|^2+\|\dbar^*\varphi \|^2+\|\varphi\|^2)
\end{equation*}
for all $\varphi \in {\mathcal D}^{0,1}_U$. Here, ${\mathcal D}^{0,1}_U$ denotes the space of $(0,1)$-forms $\varphi \in \Dom (\dbar^*)$ such that the coefficients of $\varphi$ are $C^\infty$ with compact support in $U\cap \overline \Omega$. The norm $\|\cdot\|_\varepsilon$ is the Sobolev $\varepsilon$-norm, taken coefficient-wise.
\end{definition}
A subelliptic estimate of order $1$ always holds for forms that are compactly supported in the interior of $\Omega$, i.e.\ at $p\in \Omega$. For a boundary point $p_0\in\partial \Omega$, a subelliptic estimate of order $\frac 1 2$ is known to hold if $p_0$ is a strongly pseudoconvex boundary point (\cite{Kohn_Harm_Int_I,Kohn_Harm_Int_II}). It was proven by Catlin (\cite{C1,C2,C3}) that  a subelliptic estimate holds at $p_0$ if and only if $p_0$ is of {\it finite type} in the following sense. This notion was first introduced by D'Angelo (\cite{D'Angelo_JDG}).
\begin{definition}
The {\it type} of $\partial \Omega$ at $p_0$ is defined to be
\begin{equation*}
T(\partial \Omega,p_0)=\sup _{\gamma\in{\Gamma_{p_0}}} \frac {\ord_0(r\circ\gamma)}{\ord_0(\gamma)},
\end{equation*}
where ${\Gamma_{p_0}}$ is the set of germs of non-constant local holomorphic arcs from $(\Delta,0)$ to $(\CC^n,p_0)$. We say that $p_0$ is a boundary point of {\it finite type} if $T(\partial \Omega,p_0)< \infty$. Geometrically, the type can be viewed as the maximal normalized touching order of holomorphic arcs with the boundary. It is independent of the choice of $r$.
\end{definition}
In his series of papers, Catlin was even able to establish a {\it uniform effective} subelliptic estimate: he showed that a subelliptic estimate of order at least $t^{-n^2t^{n^2}}$ holds at a boundary point $p_0$ with $T(\partial \Omega,p_0)=t$.\par
Before Catlin's work, Kohn (\cite{Kohn_Acta}) had introduced an approach to subelliptic estimates using {\it subelliptic multiplier ideals}. The idea of multiplier ideals subsequently evolved into a fundamental concept of modern algebraic geometry, although a precise understanding of the relation between the different incarnations is still lacking. Note that Catlin's work uses appropriate families of exhaustion functions and does not make use of subelliptic multipliers. The relevant definitions are as follows.
\begin{definition}
For $p\in \overline\Omega$, let $I_{p}$ be the set of all germs of complex-valued $C^\infty$ functions $F$ for which there exist a neighborhood $U$ of $p$, $\varepsilon >0$ and $C>0$ (depending on $F$) such that 
\begin{equation*}
\|F\varphi\|^2_\varepsilon \leq C(\|\dbar \varphi\|^2+\|\dbar^*\varphi \|^2+\|\varphi\|^2)
\end{equation*}
for all $\varphi \in {\mathcal D}^{0,1}_U$. Elements of $I_{p}$ are called {\it subelliptic multipliers}.
\end{definition}
Clearly, a subelliptic estimate holds at $p$ if and only if $1\in I_{p}$. The following facts about $I_p$ are established in \cite[Theorem 1.21]{Kohn_Acta}. $I_p$ is an ideal, called the {\it subelliptic multiplier ideal}. It is equal to its own real radical, meaning that, for $k\in \NN$, $|F|^k\leq |G|$ with $G\in I_p$ implies $F\in I_p$. It always contains the function $r$ and the determinant of the Levi form. Most interestingly, \cite[Theorem 1.21, (d)]{Kohn_Acta} gives a procedure of generating new multipliers from old ones: if $F_1,\ldots,F_j$ $(1\leq j\leq n-1)$ are in $I_{p}$, then the coefficients of
\begin{equation*}
\partial F_1\wedge\ldots\wedge \partial F_j\wedge\partial r\wedge\dbar r\wedge (\partial \dbar r)^{n-1-j}
\end{equation*}
are also in $I_{p}$. We call this last property $(*)$.\par
Kohn's result is that in the case of $r$ being real-analytic, the finite type condition at a boundary point $p_0$ is sufficient (and necessary) for the procedure to give $1\in I_{p_0}$.  Kohn's proof is based on property $(*)$, the real radical property, and an induction on the dimension of the zero-set of a certain increasing sequence of subideals in $I_{p_0}$.\par
The question to what extend Kohn's procedure can be used to obtain a {\it uniform effective} subelliptic estimate at a boundary point similar to (or better than) Catlin's is being investigated. When applied verbatim, the procedure was recently found not to be  always effective by means of the example
\begin{equation*}
r=\rp (z_3)+ |z_1^3+z_1z_2^b|^2+|z_2|^2
\end{equation*}
($b$ any positive integer) at the origin in $\CC^3$. The author hopes to return to this problem. 
\par
\subsection{Finite type ideals}
An interesting special (but, e.g., not necessarily geometrically convex) case of the above general situation is when $p_0=0\in \CC^{n+1}$ and $r$ (around $0$) is of the form
\begin{equation}\label{special_domain}
r=\rp (z_{n+1})+\sum_{i=1}^N |h_i(z_1,\ldots,z_{n})|^2,
\end{equation}
where $N\in\NN$, and the $h_i$ are germs of holomorphic functions vanishing at the origin in $\CC^{n}$, i.e.\ $h_i\in\frakm_{\CC^{n},0}\subset \OO_{\CC^{n},0}$.\par
It is easy to see that $T(\partial \Omega,0)<\infty$ if and only if the ideal $\fraka =(h_1,\ldots,h_N)$ is $\frakm_{\CC^{n},0}$-primary, i.e.\ there is a $q\in \NN$ such that $\frakm_{\CC^{n},0}^q\subseteq \fraka\subseteq \frakm_{\CC^{n},0}$. This means that the $h_i$ vanish jointly precisely at the origin.\par
Motivated by the above, one can make the following definition.
\begin{definition}[\cite{D'Angelo_Annals}]\label{type_definition}
Let $\fraka \subseteq \OO_{\CC^n,0}$ be an ideal of germs
of holomorphic functions at the origin of $\CC^n$ with $\frakm_{\CC^{n},0}^q\subseteq\fraka\subseteq \frakm_{\CC^{n},0}$ for some
positive integer $q$. The {\it type} of $\fraka$ is the (a priori real) number
\begin{equation*} T(\fraka):=\sup_{\gamma\in{\Gamma_0}}
\inf_{f\in \fraka }\left\{{\frac{\ord_0 \left(f\circ
\gamma\right)}{\ord_0 \left(\gamma\right)}}\right\},
\end{equation*} where ${\Gamma_0}$ is the set of all germs
of non-constant local holomorphic arcs
$\gamma:\Delta\to\CC^n$ with $\gamma (0)=0$.
\end{definition} 
Equivalently, one can write
\begin{equation}\label{2}
T(\fraka)  = 
\sup_{\gamma\in{\Gamma_0}}\left\{
\frac{\ord_0\left(\gamma^*\fraka\right )}{\ord_0\left(\gamma^*\frakm\right)}\right\},
\end{equation} where for an ideal 
$\frakb\subseteq
\calo_{\CC^n,0}$,
$\ord_0(\gamma^*\frakb)$ denotes the least order of
vanishing of the pull-backs $\gamma^*(f)$ of germs of
functions $f\in\frakb$. It is well-known (and we shall
see below) that $T(\fra)$ is actually a rational
number, and that the supremum is actually a maximum.  \par
\begin{remark}
When $r$ is of the form \eqref{special_domain}, then clearly 
\begin{equation*}
T(\partial \Omega,0)=2T(\fraka),
\end{equation*}
where $\fraka=(h_1,\ldots,h_N)$.
\end{remark}
Recently there has been some interest in understanding
geometrically the curves that compute $T(\fra)$. For
example, McNeal and N\'emethi (\cite{McN_N}) proved
that if $n = 2$, then one can specify finitely many curves $\gamma_1, \ldots, \gamma_t$ so that $T(\fra)$ is given by the maximum over the
$\gamma_\alpha$ of the ratio on the right in \eqref{2}. More on the $2$-dimensional case can be found in \cite{Favre_Jonsson}. Using ideas introduced by Teissier (\cite{Teissier1,Teissier2}), Lazarsfeld and the author (\cite{HL}) proved the analogous Theorem \ref{mthm} in all dimensions, with Proposition \ref{type_reformulation} being the key component of the proof.\label{ackn} We subsequently learned that these results were already contained in the unpublished seminar notes \cite{Lejeune_Teissier} of Lejeune-Jalabert and Teissier and therefore decided not to publish our preprint. For the convenience of the reader, and as a warm-up for the more general situation dealt with in connection with the Nullstellensatz in Section \ref{eff_nsts}, the material from \cite{HL} is reproduced in Section \ref{curve_selection}.
\begin{acknowledgement}
The author wishes to thank Brian Conrad, Jeffery McNeal, Mircea Musta\c{t}\v{a}, and, most of all, Rob Lazarsfeld for valuable discussions. Bernard Teissier kindly alerted us to, and provided a copy of, \cite{Lejeune_Teissier}.
\end{acknowledgement}
\section{Curve selection for finite type ideals}\label{curve_selection}
\subsection{Statement of the curve selection theorem}
As as matter of terminology, one says that an element
$f\in\frakb\subseteq \calo_{\CC^n,0}$ is
\textit{general} if it is a general
$\CC$-linear combination of a collection of generators
of $\frakb$.
\begin{theorem}[\cite{Lejeune_Teissier}, \cite{HL}]\label{mthm} Let $\fraka=(h_1,\ldots,h_N) \subseteq \frakm_{\CC^{n},0}$ be an $\frakm_{\CC^{n},0}$-primary ideal. Choose $n-1$ general
elements $f_1\ldots,f_{n-1}\in \fraka$ and let
\begin{equation*}
C=\Zeroes(f_1,\ldots,f_{n-1})_{{\rm red}}
\end{equation*} be the germ at $0$ of the reduced curve
arising as the common zeroes of
$f_1\ldots,f_{n-1}$. Consider the decomposition
\begin{equation*} C=C_1\cup\ldots\cup C_t
\end{equation*} of $C$ into local analytic irreducible
components,  and let
\begin{equation*}
\gamma_\alpha: \Delta\to C_\alpha
\end{equation*} be the normalization of $C_\alpha$.
Then
\begin{equation*}
T(\fraka)=\max_{\alpha=1,\ldots,t}\left\{
\frac{\ord_0\left(\gamma_\alpha^*\fraka\right)}{\ord_0\left(\gamma_\alpha^*\frakm\right)}\right\}.
\end{equation*}
\end{theorem}

\subsection{Proof of the curve selection theorem}
Fix a small neighborhood $U \subset \CC^n$ of the
origin in which the given ideal is represented by an
ideal sheaf $\fra \subseteq \OO_U$.  Let 
\begin{equation*}
\Bl_\fra(U) \lra U 
\end{equation*}
be the blowing up of $\fra$. 
Concretely, using the generators $h_1, \ldots, h_N
\in \OO(U)$ of $\fra$, $\Bl_{\fra}(U)$ can be
realized as the closure in $U \times \PP^{N-1}$ of the
graph of the meromorphic mapping $U \dra \PP^{N-1}$
defined by $h_1, \ldots, h_N$.\par

Now let $X^+\lra \Bl_\fra(U)$  be the
normalization of this blowing up, and denote by $\nu :X^+
\lra U$ the natural map.  Then
$\fraka\cdot\calo_{X^+}=\calo_{X^+}(-F)$ for some
effective Cartier divisor $F$ on $X^+$. Write
\begin{equation*} F\ =\ \sum_{i=1}^sr_iE_i
\end{equation*} for the corresponding  Weil divisor,
 where $E_i$ are irreducible divisors on $X^+$. Thus
$r_i=\ord_{E_i}(\fraka)$ is the vanishing order along
$E_i$ of the pull-back $\nu^*(f)$ of a general element
$f\in\fraka$. We also let $m_i=\ord_{E_i}(\frakm)$,
$\frakm$ denoting the maximal ideal sheaf defining the origin $0 \in U$. 
\begin{proposition}[\cite{Lejeune_Teissier}, \cite{HL}]\label{type_reformulation} With the
above notation, one has
\begin{equation}
\label{Type.Eqn} T(\fraka)=\max_{i=1,\ldots,s}\left\{
\frac{r_i}{m_i}\right\}.
\end{equation}
\end{proposition} To prove the Proposition, we will
make use of the notion of integral closure of
ideals.  The reader may consult \cite[Section
9.6.A]{PAG2} for a geometrically-oriented overview, or 
\cite{Teissier1,Teissier2} for a more detailed treatment. For our
purposes, it is sufficient to define the {\it integral
closure} of integer powers $\fraka^k$ of 
$\fraka$ to be
\begin{equation*}
\overline {\fraka^k} = \nu_*\calo_{X^+}(-kF).
\end{equation*} An important fact is the following
valuative criterion for membership in the integral
closure of an ideal (\cite[Corollaire 2, p.
328]{Teissier2}).
\begin{lemma}\label{val_criterion} Given an ideal
$\frakb\subseteq \OO_U$ vanishing only at the origin,
one has
$f\in\overline \frakb$ if and only if
\begin{equation*}
\ord_0\left(\gamma^* f\right )\geq \ord_0\left(\gamma^* \frakb\right)
\end{equation*} for all curves $\gamma\in{\Gamma_0}$. \qed
\end{lemma} 
In particular, the type of an ideal only depends on the integral closure, 
i.e.\ for $\frakb\subseteq
\calo_{\CC^n,0}$, one has $T(\frakb)= T(\overline \frakb)$. Furthermore, Lemma \ref{val_criterion} yields the
following statement, whose simple proof is left to the reader. \begin{lemma}\label{power_m_variant} The type $T(\fra)$
is the least rational number $t > 0$ such that
\begin{equation*}
\frakm^{kt}\ \subseteq\
\overline{\fraka^k}
\end{equation*}  
for all sufficiently large and divisible integers $k \in\NN$. \qed
\end{lemma}\par
 
\begin{proof}[Proof of Proposition \ref{type_reformulation}] Let $\tau$ denote the
maximum on the right hand side in
equation \eqref{Type.Eqn}. Then for all
$i=1,\ldots,s$, we have $\tau m_i\geq r_i$, which
implies that
\begin{equation*}
\left(\sum_{i=1}^s k\tau m_iE_i\right)-kF
\end{equation*} is an effective divisor on $X^+$. Let
$k$ be such that $k\tau\in
\NN$. Then
\begin{equation*}
\frakm^{k\tau}
\subseteq\nu_*\calo_{X^+}\left(-\sum_{i=1}^sk\tau
m_iE_i\right)\subseteq\nu_*\calo_{X^+}\left(-kF\right)=\overline{\fraka^k}.
\end{equation*} It follows from Lemma
\ref{power_m_variant} that $\tau\geq T(\fraka)$. On the other hand, Lemma \ref{power_m_variant} also gives
$\frakm^{kT(\fraka)}\subseteq\overline{\fraka^k}$ for
sufficiently large and divisible $k$. Therefore,
\begin{equation*} kT(\fraka)\ord_{E_i}(\frakm)\geq
k\ord_{E_i}(\fraka)
\end{equation*} for all $i=1,\ldots,s$. This implies
$\tau\leq T(\fraka)$.
\end{proof} We are now in a position to prove Theorem
\ref{mthm}.
\begin{proof}[Proof of Theorem \ref{mthm}.] Given $n-1$
general elements $f_1,\ldots,f_{n-1}\in\fraka$, we write
\begin{equation*} D_j= \divisor(f_j),\quad
\nu^*{D_j}=F+D_j',
\end{equation*} so that $D_j'$ is the proper transform
of $D_j$. Note that the
$D_j'$ generate a basepoint-free subseries of $|-F|$.
Therefore, by Bertini, the curve
\begin{equation*}
\Lambda=D_1'\cap\ldots\cap D_{n-1}'
\end{equation*} is smooth (but possibly reducible) and meets
$E=E_1\cup\ldots\cup E_s$ transversely at smooth points
of intersection of $\Lambda$ and $E$. Given an
irreducible component $C_\alpha$ of $C$, let
$\Lambda_\alpha$ denote its inverse image in $\Lambda$,
and denote by
$\gamma_\alpha:\Delta\cong\Lambda_\alpha\to C_\alpha$
the natural map. The component $\Lambda_\alpha$ is a smooth
irreducible arc that meets
$E$ transversely at one point lying on a single 
irreducible component of $E$, say
$E_{i(\alpha)}$. One also has
\begin{equation*}
\ord_0\left(\gamma^*_\alpha \fraka\right)=r_{i(\alpha)},\quad
\ord_0\left(\gamma^*_\alpha
\frakm\right)=m_{i(\alpha)}.
\end{equation*} In view of Proposition
\ref{type_reformulation}, we are done if we can show
that for each $i=1,\ldots,s$, there is an $\alpha$ such
that $i(\alpha)=i$.\par 
Now given $i$, the number of
intersections of $E_i$ with $\Lambda$ is
\begin{equation}\label{intersection_number}
\big(E_i\cdot\Lambda\big)=\big(E_i\cdot (-F)^{n-1}\big).
\end{equation} But $\calo_{X^+}(-F)$ is relatively
ample for $\nu$ (Lemma \ref{Amplitude.Lemma}), so for
all
$i$ the intersection number in
\eqref{intersection_number} is strictly positive, which
concludes the proof.
\end{proof} For the sake of completeness, we include the
key ampleness assertion from the proof as a final lemma in this section.
\begin{lemma} \label{Amplitude.Lemma}
With the above notation,
$\calo_{X^+}(-F)$ is ample on any subvariety contained
in the support of $F$.
\end{lemma}
\begin{proof} In fact, $F$ is the pull-back
under the finite mapping $X^+ \lra \Bl_\fra(U)$ of the
exceptional divisor $F_0$ on $Y = \Bl_\fra(U)$. As
explained above, we may view $Y$ as an
analytic subvariety of $U \times \PP^{r-1}$, and
$\OO_Y(-F_0)$ is the pull-back of the hyperplane line
bundle on $\PP^{r-1}$, and it is ample on $F_0$.
Since ampleness is preserved under pulling back by finite maps, the
statement follows.
\end{proof}

\section{The effective Nullstellensatz}\label{eff_nsts}
In this section, we will generalize the ideas of the previous section to obtain, in a very transparent manner, improved local and global versions of the effective Nullstellensatz, using the Theorem of Brian\c con-Skoda \cite{BS}. Ein-Lazarsfeld (\cite{Ein_Lazarsfeld}) and Hickel (\cite{Hickel}), among others, did previous work in this direction. In Subsection \ref{nsts_curve_selection}, a curve selection theorem analogous to Theorem \ref{mthm} is proven. 
\subsection{Local effective Nullstellensatz}
We now take 
$\fraka\subseteq \frakm_{\CC^{n},0}$ to be any ideal, i.e.\ we no longer require it to be $\frakm_{\CC^{n},0}$-primary. Again, we fix a small neighborhood $U \subset \CC^n$ of the origin in which the given ideal is represented by an
ideal subsheaf $\fra \subseteq \OO_U$. \par
We consider the zero-locus
\begin{equation*}
Z=\Zeroes(\fraka)_{\text{red}} \subset U
\end{equation*} 
of $\fraka$, viewed as a reduced subvariety of $U$. Let $I_Z$ be the ideal of $Z$. According to R\"uckert's Nullstellensatz, $I_Z=\sqrt{\fraka}$. Equivalently, one can say that there exists a $\sigma\in \NN$ with $I_Z^\sigma\subseteq \fraka$. We now generalize the notion of the type of an $\frakm$-primary ideal discussed above to the following invariant.
\begin{definition}
For $\fraka$ as above, we define the invariant
\begin{equation*}
\TT(\fraka)=\sup_{\gamma\in\Gamma} \left\{\frac{\ord_0\left(\gamma^*\fraka\right )}{\ord_0\left(\gamma^* I_Z\right)}\right\},
\end{equation*}
where $\Gamma$ is the set of local holomorphic arcs from $(\Delta,0)$ to $U$ whose image intersects $Z$, but is not contained in it.
\end{definition}
\begin{theorem}\label{loc_eff_nsts} $\TT(\fraka)$ governs the effective Nullstellensatz in the following way.
\begin{equation*}
\left(\sqrt{\fraka}\right)^{{\lceil n\TT(\fraka)\rceil}}= I_Z^{{\lceil n\TT(\fraka)\rceil}}\ \subseteq\ \fraka.
\end{equation*}
\end{theorem}
To prove the Theorem, we note that, in analogy to Lemma \ref{power_m_variant}, the valuative criterion gives us the following. 
\begin{lemma}\label{asympt_lemma}
$\TT(\fra)$ is the least rational number $t > 0$ such that
\begin{equation*}
I_Z^{kt}\ \subseteq\ \overline{\fraka^k}
\end{equation*}  
for all sufficiently large and divisible integers $k \in\NN$. Moreover, for any positive integer $k$, we have
\begin{equation*}
I_Z^{\lceil k\TT(\fraka)\rceil}\ \subseteq\ \overline{\fraka^k}. \qed
\end{equation*}
\end{lemma}
\begin{proof}[Proof of Theorem \ref{loc_eff_nsts}]
Apply Lemma \ref{asympt_lemma} with $k=n$ to get $I_Z^{\lceil n\TT(\fraka)\rceil}\subseteq \overline{\fraka^n}$. However, the Theorem of Brian\c con-Skoda \cite{BS} states that $\overline{\fraka^n} \subseteq \fraka$ and the Theorem is proven.
\end{proof}
\subsection{A curve selection theorem}\label{nsts_curve_selection}
In the present situation, we can reintroduce the notation from the previous section: Let $\Bl_\fra(U) \lra U$ be the blowing up of $\fra$ and $X^+\lra \Bl_\fra(U)$  the normalization of this blowing up. Denote by $\nu :X^+
\lra U$ the natural map.  Then again
$\fraka\cdot\calo_{X^+}=\calo_{X^+}(-F)$ for some
effective Cartier divisor $F$ on $X^+$. Again, write $F\ =\ \sum_{i=1}^sr_iE_i$ for the corresponding  Weil divisor, where $E_i$ are irreducible divisors on $X^+$. Thus
$r_i=\ord_{E_i}(\fraka)$ is the vanishing order along
$E_i$ of the pull-back $\nu^*(f)$ of a general element
$f\in\fraka$. We also let $m_i=\ord_{E_i}(I_Z)$. In the same way in which Lemma \ref{power_m_variant} was used to prove Proposition \ref{type_reformulation}, Lemma \ref{asympt_lemma} can be used to establish the following Proposition.
\begin{proposition}\label{tt_as_ratio}
\begin{equation}
\TT(\fraka)=\max_{i=1,\ldots,s}\left\{\frac{r_i}{m_i}\right\}. \qed
\end{equation}
\end{proposition}
We would like to describe a finite set of curves that compute $\TT(\fraka)$. This is done in the following Theorem. 
\begin{theorem}
Let $d$ denote the dimension of $Z$ (i.e.\ the maximal dimension of its irreducible components).
Take $f_1,\ldots,f_{n-1}\in \fraka$ and $g_1,\ldots,g_{d}\in \OO_U$ to be general elements. Let 
\begin{eqnarray*}
C^{0}&=&\overline{\Zeroes(f_1,\ldots,f_{n-1})_{{\rm red}}\cap (U\backslash Z)}\\
C^{1}&=&\overline{\Zeroes(f_1,\ldots,f_{n-2},g_1)_{{\rm red}}\cap (U\backslash Z)}\\
&\vdots&\\
C^{d}&=&\overline{\Zeroes(f_1,\ldots,f_{n-1-d},g_1\ldots,g_{d})_{{\rm red}}\cap (U\backslash Z)}
\end{eqnarray*}
be $d$ curves (possibly reducible or empty). For $\delta=0,\ldots,d$, let $p^\delta_1,\ldots,p^\delta_{M_\delta}$ be the points of intersection of $C^{\delta}$ with $Z$. In small neighborhoods $U(p^\delta_{\mu_\delta})$ of each $p^\delta_{\mu_\delta}$ ($\mu_\delta=1,\ldots,M_\delta$), we take an irreducible decomposition
\begin{equation*}
U(p^\delta_{\mu_\delta})\cap C^\delta=C^{\delta,\mu_\delta}_1\cup\ldots\cup C^{\delta,\mu_\delta}_{t_{\delta,\mu_\delta}}.
\end{equation*}
For $\alpha_{\delta,\mu_\delta}=1,\ldots,t_{\delta,\mu_\delta}$, we take normalization maps 
\begin{equation*}
\gamma^{\delta,\mu_\delta}_{\alpha_{\delta,\mu_\delta}}:\Delta\to C^{\delta,\mu_\delta}_{\alpha_{\delta,\mu_\delta}}.
\end{equation*}
Then 
\begin{equation*}
\TT(\fraka)=\max_{\delta=0\ldots,d;\,\mu_\delta=1,\ldots,M_\delta;\,\alpha_{\delta,\mu_\delta}=1,\ldots,t_{\delta,\mu_\delta}} \left\{\frac{\ord_0\bigl(\bigl(\gamma^{\delta,\mu_\delta}_{\alpha_{\delta,\mu_\delta}}\bigr)^*\fraka\bigr)}{\ord_0\bigl(\bigl(\gamma^{\delta,\mu_\delta}_{\alpha_{\delta,\mu_\delta}}\bigr)^* I_Z\bigr)}\right\}.
\end{equation*}
\end{theorem}
\begin{proof}
For $j=1,\ldots,n-1$, write
\begin{equation*}
D_j=\divisor(f_j),\quad \nu^*D_j=F+D_j',
\end{equation*}
and for $\iota=1,\ldots,d$
\begin{equation*}
G_\iota=\divisor(g_\iota),\quad \nu^*G_\iota=G_\iota'.
\end{equation*}
For $\delta=0,\ldots,d$, let
\begin{equation*}
\Lambda^\delta=D_1'\cap\ldots\cap D_{n-1-\delta}'\cap G_1'\cap\ldots\cap G_\delta'.
\end{equation*}
By Bertini's Theorem, $\Lambda^\delta$ is a smooth (but possibly reducible) curve which meets $E=E_1\cup\ldots\cup E_s$ transversely at smooth points of intersection with $E$. Given an irreducible component $C^{\delta,\mu_\delta}_{\alpha_{\delta,\mu_\delta}}$, let $\Lambda^{\delta,\mu_\delta}_{\alpha_{\delta,\mu_\delta}}$ be the component of $\Lambda^\delta$ that maps to it. We denote the natural isomorphism by $\gamma^{\delta,\mu_\delta}_{\alpha_{\delta,\mu_\delta}}:\Delta\cong\Lambda^{\delta,\mu_\delta}_{\alpha_{\delta,\mu_\delta}}\to C^{\delta,\mu_\delta}_{\alpha_{\delta,\mu_\delta}}$. \par
If $\Lambda^{\delta,\mu_\delta}_{\alpha_{\delta,\mu_\delta}}$ meets $E$ in the component $E_{i(\alpha_{\delta,\mu_\delta})}$, then
\begin{equation*}
\ord_0\bigl(\bigl(\gamma^{\delta,\mu_\delta}_{\alpha_{\delta,\mu_\delta}}\bigr)^*\fraka\bigr)=r_{i(\alpha_{\delta,\mu_\delta})},\quad \ord_0\bigl(\bigl(\gamma^{\delta,\mu_\delta}_{\alpha_{\delta,\mu_\delta}}\bigr)^*I_Z\bigr)=m_{i(\alpha_{\delta,\mu_\delta})}.
\end{equation*}
From now on, we work with an arbitrary but fixed $\delta$. In view of Proposition \ref{tt_as_ratio}, it suffices to show that for all $E_i$ with $\dim \nu(E_i)=\delta$, there exists a $\Lambda^{\delta,\mu_\delta}_{\alpha_{\delta,\mu_\delta}}$ which meets $E_i$.\par
The intersection of $\nu(E_i)$ with $g_1\ldots,g_\delta$ is a (nonempty) set of isolated points, and we can concentrate on any one of them, say $p$. The intersection number of $E_i\cap {\nu^{-1}(p)}$ with $\Lambda^\delta$ is
\begin{equation}\label{inters_with_generic_fiber}
(E_i \cdot {\nu^{-1}(p)} \cdot (-F)^{n-1-\delta}).
\end{equation}
In the same way as in Lemma \ref{Amplitude.Lemma}, one can see that $\calo_{X^+}(-F)|_{E_i\cap\nu^{-1}(p)}$ is ample, so the intersection number in \eqref{inters_with_generic_fiber} is strictly positive. 
\end{proof}
\begin{remark}
Based merely on the definition of the curves $C^\delta$, it is not obvious that they are always nonempty when they are needed to be nonempty. Note that the positivity of \eqref{inters_with_generic_fiber} proves this fact in passing.
\end{remark}

\subsection{Geometric effective Nullstellensatz}\label{global_geom_nsts}
Let $X$ be a non-singular projective variety. In this subsection, we take $\fraka$ to be an ideal subsheaf of $\OO_X$. Let $L$ be an ample divisor on $X$ such that $\OO_X(L)\otimes\fraka$ is globally generated.\par
We extend the definitions of $Z$ and $\TT(\fraka)$ in the obvious way to the the present global situation, simply by replacing $U$ by $X$. Note that the inclusion $I_Z^{{\lceil n\TT(\fraka)\rceil}}\ \subseteq\ \fraka$ established in Theorem \ref{loc_eff_nsts} continues to hold true (because it is verified locally). The purpose of this subsection is to establish the fact that $\TT(\fraka)$ is bounded by the top self-intersection of $L$:
\begin{theorem} \label{tt_less_than_Ln} With the above notation,
\begin{equation*}
\TT(\fraka)\leq (L^n).
\end{equation*}
\end{theorem}
We immediately obtain
\begin{corollary}[{\cite[p. 431]{Ein_Lazarsfeld}}]
With the above notation, 
\begin{equation*}
\left(\sqrt{\fraka}\right)^{n(L^n)}=I_Z^{n(L^n)}\subseteq \fraka.
\end{equation*}
\end{corollary}
\begin{proof}[Proof of Theorem \ref{tt_less_than_Ln}]
Let $\Bl_\fra(X) \lra X$ be the blowing up of $\fra$ and $X^+\lra \Bl_\fra(X)$  the
normalization of this blowing up. Denote by $\nu :X^+
\lra X$ the natural map.  Then again
$\fraka\cdot\calo_{X^+}=\calo_{X^+}(-F)$ for some
effective Cartier divisor $F$ on $X^+$. Again, write $F\ =\ \sum_{i=1}^sr_iE_i$ for the corresponding Weil divisor, where $E_i$ are irreducible divisors on $X^+$. Thus
$r_i=\ord_{E_i}(\fraka)$ is the vanishing order along
$E_i$ of the pull-back $\nu^*(f)$ of a general element
$f\in\fraka$. We also let $m_i=\ord_{E_i}(I_Z)$. For the same reason as before, $\TT(\fraka)=\max_{i=1,\ldots,s}\left\{\frac{r_i}{m_i}\right\}$.\par
The following estimate is based on the telescoping sum trick used in \cite[Proposition 3.1]{Ein_Lazarsfeld}. Let $M$ be a divisor on $X^+$ linearly equivalent to $\nu^*L-F$. Note that $M$ is nef. We have
\begin{eqnarray*}
(L^n)&=&((\nu^* L)^n)\\
&\geq&((\nu^* L)^n)-(M^n)\\
&=&\bigl(\bigl(\nu^* L-M\bigr)\cdot\bigl(\sum_{\mu=0}^{n-1} (\nu^* L)^{\mu}\cdot M^{n-1-\mu}\bigr) \bigr)\\
&=&\bigl(F\cdot \bigl(\sum_{\mu=0}^{n-1} (\nu^* L)^{\mu}\cdot M^{n-1-\mu}\bigr)\bigr)\\
&=&\sum_{i=1}^s\bigl(r_iE_i \cdot \bigl(\sum_{\mu=0}^{n-1} (\nu^* L)^{\mu}\cdot M^{n-1-\mu}\bigr)\bigr)\\
&\geq&\sum_{i=1}^s\bigl(r_iE_i \cdot (\nu^* L)^{\dim \nu(E_i) }\cdot M^{n-1-\dim \nu(E_i)}\bigr)\\
&=&\sum_{i=1}^s(r_iE_i \cdot (\nu^* L)^{\dim \nu(E_i) }\cdot (\nu^* L-F)^{n-1-\dim \nu(E_i)})\\
&\geq& \sum_{i=1}^s(r_iE_i \cdot (\nu^* L)^{\dim \nu(E_i) }\cdot (-F)^{n-1-\dim \nu(E_i)})\\
&\geq& \sum_{i=1}^s r_i.
\end{eqnarray*}
Note that for the last inequality we are again exploiting the $\nu$-ampleness of $-F$, which makes the intersection number that is the coefficient of $r_i$ strictly positive for all $i$. On the whole, we have
\begin{equation*}
\TT(\fraka)=\max_{i=1,\ldots,s}\left\{\frac{r_i}{m_i}\right\}\leq  \sum_{i=1}^s r_i\leq (L^n).\qedhere
\end{equation*}
\end{proof}

\pagebreak


\begin{thebibliography}{Laz04b}

\bibitem[Cat83]{C1}
D.~Catlin.
\newblock Necessary conditions for subellipticity of the {$\bar \partial
  $}-{N}eumann problem.
\newblock {\em Ann. of Math. (2)}, 117(1):147--171, 1983.

\bibitem[Cat84]{C2}
D.~Catlin.
\newblock Boundary invariants of pseudoconvex domains.
\newblock {\em Ann. of Math. (2)}, 120(3):529--586, 1984.

\bibitem[Cat87]{C3}
D.~Catlin.
\newblock Subelliptic estimates for the {$\overline\partial$}-{N}eumann problem
  on pseudoconvex domains.
\newblock {\em Ann. of Math. (2)}, 126(1):131--191, 1987.

\bibitem[D'A79]{D'Angelo_JDG}
J.~P. D'Angelo.
\newblock Finite type conditions for real hypersurfaces.
\newblock {\em J. Differential Geom.}, 14(1):59--66 (1980), 1979.

\bibitem[D'A82]{D'Angelo_Annals}
J.~P. D'Angelo.
\newblock Real hypersurfaces, orders of contact, and applications.
\newblock {\em Ann. of Math. (2)}, 115(3):615--637, 1982.

\bibitem[D'A93]{D'Angelo-book}
J.~P. D'Angelo.
\newblock {\em Several complex variables and the geometry of real
  hypersurfaces}.
\newblock Studies in Advanced Mathematics. CRC Press, Boca Raton, FL, 1993.

\bibitem[D'A95]{DAngelo_Gunning_Kohn_Vol}
J.~P. D'Angelo.
\newblock Finite type conditions and subelliptic estimates.
\newblock In {\em Modern methods in complex analysis (Princeton, NJ, 1992)},
  volume 137 of {\em Ann. of Math. Stud.}, pages 63--78. Princeton Univ. Press,
  Princeton, NJ, 1995.

\bibitem[DK99]{D'Angelo-Kohn}
J.~P. D'Angelo and J.~J. Kohn.
\newblock Subelliptic estimates and finite type.
\newblock In {\em Several complex variables (Berkeley, CA, 1995--1996)},
  volume~37 of {\em Math. Sci. Res. Inst. Publ.}, pages 199--232. Cambridge
  Univ. Press, Cambridge, 1999.

\bibitem[EL99]{Ein_Lazarsfeld}
L.~Ein and R.~Lazarsfeld.
\newblock A geometric effective {N}ullstellensatz.
\newblock {\em Invent. Math.}, 137(2):427--448, 1999.

\bibitem[FJ04]{Favre_Jonsson}
Ch. Favre and M.~Jonsson.
\newblock {\em The valuative tree}, volume 1853 of {\em Lecture Notes in
  Mathematics}.
\newblock Springer-Verlag, Berlin, 2004.

\bibitem[Hic01]{Hickel}
M.~Hickel.
\newblock Solution d'une conjecture de {C}.\ {B}erenstein--{A}.\ {Y}ger et
  invariants de contact \`a l'infini.
\newblock {\em Ann. Inst. Fourier (Grenoble)}, 51(3):707--744, 2001.

\bibitem[HL05]{HL}
G.~Heier and R.~Lazarsfeld.
\newblock Curve selection for finite-type ideals.
\newblock {\em {\rm math.CV/0506557}}, 2005.

\bibitem[KN65]{Kohn_Nirenberg}
J.~J. Kohn and L.~Nirenberg.
\newblock Non-coercive boundary value problems.
\newblock {\em Comm. Pure Appl. Math.}, 18:443--492, 1965.

\bibitem[Koh63]{Kohn_Harm_Int_I}
J.~J. Kohn.
\newblock Harmonic integrals on strongly pseudo-convex manifolds. {I}.
\newblock {\em Ann. of Math. (2)}, 78:112--148, 1963.

\bibitem[Koh64]{Kohn_Harm_Int_II}
J.~J. Kohn.
\newblock Harmonic integrals on strongly pseudo-convex manifolds. {II}.
\newblock {\em Ann. of Math. (2)}, 79:450--472, 1964.

\bibitem[Koh79]{Kohn_Acta}
J.~J. Kohn.
\newblock Subellipticity of the {$\bar \partial $}-{N}eumann problem on
  pseudo-convex domains: sufficient conditions.
\newblock {\em Acta Math.}, 142(1-2):79--122, 1979.

\bibitem[Laz04a]{PAG2}
R.~Lazarsfeld.
\newblock {\em Positivity in algebraic geometry. {I}}, volume~48 of {\em
  Ergebnisse der Mathematik und ihrer Grenzgebiete. 3. Folge. A Series of
  Modern Surveys in Mathematics [Results in Mathematics and Related Areas. 3rd
  Series. A Series of Modern Surveys in Mathematics]}.
\newblock Springer-Verlag, Berlin, 2004.
\newblock Classical setting: line bundles and linear series.

\bibitem[Laz04b]{PAG1}
R.~Lazarsfeld.
\newblock {\em Positivity in algebraic geometry. {II}}, volume~49 of {\em
  Ergebnisse der Mathematik und ihrer Grenzgebiete. 3. Folge. A Series of
  Modern Surveys in Mathematics [Results in Mathematics and Related Areas. 3rd
  Series. A Series of Modern Surveys in Mathematics]}.
\newblock Springer-Verlag, Berlin, 2004.
\newblock Positivity for vector bundles, and multiplier ideals.

\bibitem[LJT74]{Lejeune_Teissier}
M.~Lejeune-Jalabert and B.~Teissier.
\newblock Cl\^{o}ture int\'egral des id\'eaux et \'equi\-singularit\'e.
\newblock {\em Grenoble seminar notes (to appear in Ann. Fac. Sci. Toulouse
  Math. (6))}, 1974.

\bibitem[MN05]{McN_N}
J.~D. McNeal and A.~N{\'e}methi.
\newblock The order of contact of a holomorphic ideal in {${\mathbb {C}}^2$}.
\newblock {\em Math. Z.}, 250(4):873--883, 2005.

\bibitem[SB74]{BS}
H.~Skoda and J.~Brian{\c{c}}on.
\newblock Sur la cl\^oture int\'egrale d'un id\'eal de germes de fonctions
  holomorphes en un point de {${\bf C}\sp{n}$}.
\newblock {\em C. R. Acad. Sci. Paris S\'er. A}, 278:949--951, 1974.

\bibitem[Siu05]{Siu}
Y.-T. Siu.
\newblock Multiplier ideal sheaves in complex and algebraic geometry.
\newblock {\em Sci. China Ser. A}, 48(suppl.):1--31, 2005.

\bibitem[Tei77]{Teissier1}
B.~Teissier.
\newblock Vari\'et\'es polaires. {I}. {I}nvariants polaires des singularit\'es
  d'hyper\-surfaces.
\newblock {\em Invent. Math.}, 40(3):267--292, 1977.

\bibitem[Tei82]{Teissier2}
B.~Teissier.
\newblock Vari\'et\'es polaires. {II}. {M}ultiplicit\'es polaires, sections
  planes, et conditions de {W}hitney.
\newblock In {\em Algebraic geometry (La R\'abida, 1981)}, volume 961 of {\em
  Lecture Notes in Math.}, pages 314--491. Springer, Berlin, 1982.

\end{thebibliography}
\end{document}